\documentclass[12pt]{article}

\usepackage{amsmath,amssymb,amsthm,amscd,a4wide,upref}

\usepackage[enableskew]{youngtab}

\usepackage{hyperref}
\hypersetup{colorlinks,citecolor=blue,filecolor=black,linkcolor=blue,urlcolor=blue}
\usepackage{enumerate}
\usepackage{mathtools}
\usepackage{enumitem}

\usepackage{lmodern}     
\usepackage[T1]{fontenc}

\begin{document}


\newcommand{\ad}{{\rm ad}}
\newcommand{\cri}{{\rm cri}}
\newcommand{\row}{{\rm row}}
\newcommand{\col}{{\rm col}}
\newcommand{\Ann}{{\rm{Ann}\ts}}
\newcommand{\End}{{\rm{End}\ts}}
\newcommand{\Rep}{{\rm{Rep}\ts}}
\newcommand{\Hom}{{\rm{Hom}}}
\newcommand{\Mat}{{\rm{Mat}}}
\newcommand{\ch}{{\rm{ch}\ts}}
\newcommand{\chara}{{\rm{char}\ts}}
\newcommand{\diag}{{\rm diag}}
\newcommand{\st}{{\rm st}}
\newcommand{\non}{\nonumber}
\newcommand{\wt}{\widetilde}
\newcommand{\wh}{\widehat}
\newcommand{\ol}{\overline}
\newcommand{\ot}{\otimes}
\newcommand{\la}{\lambda}
\newcommand{\La}{\Lambda}
\newcommand{\De}{\Delta}
\newcommand{\al}{\alpha}
\newcommand{\be}{\beta}
\newcommand{\ga}{\gamma}
\newcommand{\Ga}{\Gamma}
\newcommand{\ep}{\epsilon}
\newcommand{\ka}{\kappa}
\newcommand{\vk}{\varkappa}
\newcommand{\si}{\sigma}
\newcommand{\vs}{\varsigma}
\newcommand{\vp}{\varphi}
\newcommand{\ta}{\theta}
\newcommand{\de}{\delta}
\newcommand{\ze}{\zeta}
\newcommand{\om}{\omega}
\newcommand{\Om}{\Omega}
\newcommand{\ee}{\epsilon^{}}
\newcommand{\su}{s^{}}
\newcommand{\hra}{\hookrightarrow}
\newcommand{\ve}{\varepsilon}
\newcommand{\pr}{^{\tss\prime}}
\newcommand{\ts}{\,}
\newcommand{\vac}{\mathbf{1}}
\newcommand{\vacu}{|0\rangle}
\newcommand{\di}{\partial}
\newcommand{\qin}{q^{-1}}
\newcommand{\tss}{\hspace{1pt}}
\newcommand{\Sr}{ {\rm S}}
\newcommand{\U}{ {\rm U}}
\newcommand{\BL}{ {\overline L}}
\newcommand{\BE}{ {\overline E}}
\newcommand{\BP}{ {\overline P}}
\newcommand{\AAb}{\mathbb{A}\tss}
\newcommand{\CC}{\mathbb{C}\tss}
\newcommand{\KK}{\mathbb{K}\tss}
\newcommand{\QQ}{\mathbb{Q}\tss}
\newcommand{\SSb}{\mathbb{S}\tss}
\newcommand{\TT}{\mathbb{T}\tss}
\newcommand{\ZZ}{\mathbb{Z}\tss}
\newcommand{\DY}{ {\rm DY}}
\newcommand{\X}{ {\rm X}}
\newcommand{\Y}{ {\rm Y}}
\newcommand{\Z}{{\rm Z}}
\newcommand{\ZX}{{\rm ZX}}
\newcommand{\ZY}{{\rm ZY}}
\newcommand{\Ac}{\mathcal{A}}
\newcommand{\Lc}{\mathcal{L}}
\newcommand{\Mc}{\mathcal{M}}
\newcommand{\Pc}{\mathcal{P}}
\newcommand{\Qc}{\mathcal{Q}}
\newcommand{\Rc}{\mathcal{R}}
\newcommand{\Sc}{\mathcal{S}}
\newcommand{\Tc}{\mathcal{T}}
\newcommand{\Bc}{\mathcal{B}}
\newcommand{\Ec}{\mathcal{E}}
\newcommand{\Fc}{\mathcal{F}}
\newcommand{\Gc}{\mathcal{G}}
\newcommand{\Hc}{\mathcal{H}}
\newcommand{\Uc}{\mathcal{U}}
\newcommand{\Vc}{\mathcal{V}}
\newcommand{\Wc}{\mathcal{W}}
\newcommand{\Yc}{\mathcal{Y}}
\newcommand{\Cl}{\mathcal{C}l}
\newcommand{\Ar}{{\rm A}}
\newcommand{\Br}{{\rm B}}
\newcommand{\Ir}{{\rm I}}
\newcommand{\Fr}{{\rm F}}
\newcommand{\Jr}{{\rm J}}
\newcommand{\Or}{{\rm O}}
\newcommand{\GL}{{\rm GL}}
\newcommand{\Spr}{{\rm Sp}}
\newcommand{\Rr}{{\rm R}}
\newcommand{\Zr}{{\rm Z}}
\newcommand{\gl}{\mathfrak{gl}}
\newcommand{\middd}{{\rm mid}}
\newcommand{\ev}{{\rm ev}}
\newcommand{\Pf}{{\rm Pf}}
\newcommand{\Norm}{{\rm Norm\tss}}
\newcommand{\oa}{\mathfrak{o}}
\newcommand{\spa}{\mathfrak{sp}}
\newcommand{\osp}{\mathfrak{osp}}
\newcommand{\f}{\mathfrak{f}}
\newcommand{\se}{\mathfrak{s}}
\newcommand{\g}{\mathfrak{g}}
\newcommand{\h}{\mathfrak h}
\newcommand{\n}{\mathfrak n}
\newcommand{\m}{\mathfrak m}
\newcommand{\z}{\mathfrak{z}}
\newcommand{\Zgot}{\mathfrak{Z}}
\newcommand{\p}{\mathfrak{p}}
\newcommand{\sll}{\mathfrak{sl}}
\newcommand{\agot}{\mathfrak{a}}
\newcommand{\bgot}{\mathfrak{b}}
\newcommand{\qdet}{ {\rm qdet}\ts}
\newcommand{\Ber}{ {\rm Ber}\ts}
\newcommand{\HC}{ {\mathcal HC}}
\newcommand{\cdet}{{\rm cdet}}
\newcommand{\rdet}{{\rm rdet}}
\newcommand{\tr}{ {\rm tr}}
\newcommand{\gr}{ {\rm gr}\ts}
\newcommand{\str}{ {\rm str}}
\newcommand{\loc}{{\rm loc}}
\newcommand{\Gr}{{\rm G}}
\newcommand{\sgn}{ {\rm sgn}\ts}
\newcommand{\sign}{{\rm sgn}}
\newcommand{\ba}{\bar{a}}
\newcommand{\bb}{\bar{b}}
\newcommand{\bi}{\bar{\imath}}
\newcommand{\bj}{\bar{\jmath}}
\newcommand{\bk}{\bar{k}}
\newcommand{\bl}{\bar{l}}
\newcommand{\bp}{\bar{p}}
\newcommand{\hb}{\mathbf{h}}
\newcommand{\Sym}{\mathfrak S}
\newcommand{\fand}{\quad\text{and}\quad}
\newcommand{\Fand}{\qquad\text{and}\qquad}
\newcommand{\For}{\qquad\text{or}\qquad}
\newcommand{\for}{\quad\text{or}\quad}
\newcommand{\grpr}{{\rm gr}^{\tss\prime}\ts}
\newcommand{\degpr}{{\rm deg}^{\tss\prime}\tss}
\newcommand{\bideg}{{\rm bideg}\ts}

\renewcommand{\theequation}{\arabic{section}.\arabic{equation}}

\numberwithin{equation}{section}

\newtheorem{thm}{Theorem}[section]
\newtheorem{lem}[thm]{Lemma}
\newtheorem{prop}[thm]{Proposition}
\newtheorem{cor}[thm]{Corollary}
\newtheorem{conj}[thm]{Conjecture}
\newtheorem*{mthm}{Main Theorem}
\newtheorem*{mthma}{Theorem A}
\newtheorem*{mthmb}{Theorem B}
\newtheorem*{mthmc}{Theorem C}
\newtheorem*{mthmd}{Theorem D}

\theoremstyle{definition}
\newtheorem{defin}[thm]{Definition}

\theoremstyle{remark}
\newtheorem{remark}[thm]{Remark}
\newtheorem{example}[thm]{Example}
\newtheorem{examples}[thm]{Examples}

\newcommand{\bth}{\begin{thm}}
\renewcommand{\eth}{\end{thm}}
\newcommand{\bpr}{\begin{prop}}
\newcommand{\epr}{\end{prop}}
\newcommand{\ble}{\begin{lem}}
\newcommand{\ele}{\end{lem}}
\newcommand{\bco}{\begin{cor}}
\newcommand{\eco}{\end{cor}}
\newcommand{\bde}{\begin{defin}}
\newcommand{\ede}{\end{defin}}
\newcommand{\bex}{\begin{example}}
\newcommand{\eex}{\end{example}}
\newcommand{\bes}{\begin{examples}}
\newcommand{\ees}{\end{examples}}
\newcommand{\bre}{\begin{remark}}
\newcommand{\ere}{\end{remark}}
\newcommand{\bcj}{\begin{conj}}
\newcommand{\ecj}{\end{conj}}

\newcommand{\bal}{\begin{aligned}}
\newcommand{\eal}{\end{aligned}}
\newcommand{\beq}{\begin{equation}}
\newcommand{\eeq}{\end{equation}}
\newcommand{\ben}{\begin{equation*}}
\newcommand{\een}{\end{equation*}}

\newcommand{\bpf}{\begin{proof}}
\newcommand{\epf}{\end{proof}}

\def\beql#1{\begin{equation}\label{#1}}

\newcommand{\Res}{\mathop{\mathrm{Res}}}

\title{\Large\bf Odd reflections in the Yangian associated with $\gl(m|n)$}

\author{A. I. Molev}

\date{} 
\maketitle


\begin{abstract}
The odd reflections are an effective tool in the Lie superalgebra representation
theory, as they relate non-conjugate Borel subalgebras.
We introduce analogues of the
odd reflections for the Yangian $\Y(\gl_{m|n})$ and use them to produce
a transition rule for
the parameters of the highest weight modules corresponding to a change of the
parity sequence. This leads to a description of the finite-dimensional irreducible
representations of the Yangian associated with an arbitrary
parity sequence.
\end{abstract}



%

\section{Introduction}\label{sec:int}
\setcounter{equation}{0}

The {\em Yangian} $\Y(\gl_{m|n})$ associated with the Lie superalgebra $\gl_{m|n}$ is a
deformation of the universal enveloping algebra $\U(\gl_{m|n}[u])$ in the class of
Hopf algebras. The original definition is due to Nazarov~\cite{n:qb}, where the Yangian was
introduced via an $R$-matrix presentation. A Drinfeld-type presentation corresponding to a standard Borel
subalgebra was obtained by Gow~\cite{g:gd}, extending the results of Drinfeld~\cite{d:nr}
and Brundan and Kleshchev~\cite{bk:pp} on the Yangian $\Y(\gl_n)$. An earlier note
by Stukopin~\cite{s:yl} outlines an alternative approach to
the Yangian presentations (with a correction in the set of defining relations
pointed out in \cite{g:gd}). More recently, Drinfeld-type
parabolic presentations of the Yangian
$\Y(\gl_{m|n})$, corresponding to arbitrary Borel
subalgebras in $\gl_{m|n}$ were given by Peng~\cite{p:pp}. In a particular case,
such a presentation was reproduced by
Tsymbaliuk~\cite{t:sa} who also gave a description of the center of
the Yangian by using the {\em quantum Berezinian},
thus generalizing the results of
\cite{n:qb} and \cite{g:rl}.

It was shown by Zhang~\cite{zh:sy}, that the finite-dimensional irreducible representations of
the Yangian $\Y(\gl_{m|n})$ are described in a way similar to the
representations of the Yangians
associated with the simple Lie algebras, as given by Tarasov~\cite{t:im} and
Drinfeld~\cite{d:nr}. The classification in \cite{zh:sy} relies on the presentation of the Yangian
corresponding to the standard Borel subalgebra in $\gl_{m|n}$
(or the standard parity sequence) and leaves open the question as to
how the finite-dimensionality conditions on irreducible highest weight representations
can be stated for an arbitrary Borel subalgebra.
For the Lie superalgebra $\gl_{m|n}$ itself (along with other basic Lie superalgebras),
the same question can be answered with the use of the {\em odd reflections}; see
e.g. \cite[Secs~1.3 \& 2.4]{cw:dr} and \cite[Sec.~3.5]{m:ls}. They originate in
the work of Serganova~\cite[Appendix]{lss:eo}; some quantum versions (as {\em super-reflections})
were considered by Yamane~\cite{y:dr}.

Our goal in this paper is to introduce analogues of the odd reflections for the
Yangian $\Y(\gl_{m|n})$. More precisely, we derive a transition rule
for the parameters of the irreducible highest weight representations
corresponding to
a change of the
parity sequence. Accordingly, necessary and sufficient conditions for
such representations to be finite-dimensional can be obtained by applying
a chain of transitions. As this relies on such conditions for the standard parity sequence,
we give a brief review (in Section~\ref{sec:hw}) of the
classification results for the finite-dimensional
irreducible representations of
the Yangian $\Y(\gl_{m|n})$ following \cite{zh:rs} and \cite{zh:sy}; cf. \cite{s:ry}
where an approach based on
the Drinfeld presentation was used.

The transition rule is associated with an odd simple isotropic root for a given
Borel subalgebra of $\gl_{m|n}$. This determines a subalgebra of the Yangian $\Y(\gl_{m|n})$
isomorphic to $\Y(\gl_{1|1})$. Therefore, the arguments are essentially
reduced to the analysis of
the irreducible highest weight representations of the Yangian $\Y(\gl_{1|1})$.
The key step is to describe how the odd reflections affect the highest weights.
The calculations rely on the properties of the quantum Berezinian
and they lead to the general transition rules as given in Theorem~\ref{thm:oddre}
and Corollary~\ref{cor:oddrege}.

The version of the transition rule stated in Corollary~\ref{cor:oddrege}
should apply in a similar form to the orthosymplectic Yangians introduced in \cite{aacfr:rp}.
However, apart from the particular case of $\osp_{1|2n}$ considered in
\cite{m:ry},
the classification problem for
their finite-dimensional irreducible representations is still open.

After the first version of the paper was posted in the arXiv,
one more proof of Proposition~\ref{prop:prove} was given in \cite{l:no}.

\section{Definition and basic properties of the Yangian}
\label{sec:db}

For given nonnegative integers $m$ and $n$ consider
the {\em parity sequences} $\se=\se_1\dots\se_{m+n}$ of length $m+n$,
where each term $\se_i$ is $0$ or $1$, and the total number of zeros is $m$.
The {\em standard sequence} $\se^{\st}=0\dots 0\tss 1\dots 1$
is defined by $\se_i=0$ for $i=1,\dots,m$
and $\se_i=1$ for $i=m+1,\dots,m+n$.

When a parity sequence $\se$ is fixed, we will simply write $\bi$ to
denote its $i$-th term $\se_i$. For a fixed sequence
introduce the $\ZZ_2$-graded
vector space $\CC^{m|n}$ over $\CC$ with the basis
$e_1,e_2,\dots,e_{m+n}$, where the parity of the basis vector
$e_i$ is defined to be $\bi\mod 2$.
Accordingly, equip
the endomorphism algebra $\End\CC^{m|n}$ with the $\ZZ_2$-gradation, where
the parity of the matrix unit $e_{ij}$ is found by
$\bi+\bj\mod 2$.

A standard basis of the general linear Lie superalgebra $\gl_{m|n}$ is formed by elements $E_{ij}$
of the parity $\bi+\bj\mod 2$ for $1\leqslant i,j\leqslant m+n$ with the commutation relations
\ben
[E_{ij},E_{kl}]
=\de_{kj}\ts E_{i\tss l}-\de_{i\tss l}\ts E_{kj}(-1)^{(\bi+\bj)(\bk+\bl)}.
\een
The {\em Yang $R$-matrix} associated with $\gl_{m|n}$ is the
rational function in $u$ given by
$
R(u)=1-P\tss u^{-1},
$
where $P$ is the permutation operator,
\ben
P=\sum_{i,j=1}^{m+n} e_{ij}\ot e_{ji}(-1)^{\bj}\in \End\CC^{m|n}\ot\End\CC^{m|n}.
\een
Following \cite{n:qb},
define the {\em Yangian}
$\Y(\gl^{\tss\se}_{m|n})$ associated with $\se$
(written simply as $\Y(\gl_{m|n})$, if the sequence $\se$ is fixed),
as the $\ZZ_2$-graded algebra with generators
$t_{ij}^{(r)}$ of parity $\bi+\bj\mod 2$, where $1\leqslant i,j\leqslant m+n$ and $r=1,2,\dots$,
satisfying the following quadratic relations. To write them down,
introduce the formal series
\beql{tiju}
t_{ij}(u)=\de_{ij}+\sum_{r=1}^{\infty}t_{ij}^{(r)}\ts u^{-r}
\in\Y(\gl_{m|n})[[u^{-1}]]
\eeq
and combine them into the matrix $T(u)=[t_{ij}(u)]$ so that
\ben
T(u)=\sum_{i,j=1}^{m+n} e_{ij}\ot t_{ij}(u)(-1)^{\bi\tss\bj+\bj}
\in \End\CC^{m|n}\ot \Y(\gl_{m|n})[[u^{-1}]].
\een
Consider the algebra
$\End\CC^{m|n}\ot\End\CC^{m|n}\ot \Y(\gl_{m|n})[[u^{-1}]]$
and introduce its elements $T_1(u)$ and $T_2(u)$ by
\ben
T_1(u)=\sum_{i,j=1}^{m+n} e_{ij}\ot 1\ot t_{ij}(u)(-1)^{\bi\tss\bj+\bj},\qquad
T_2(u)=\sum_{i,j=1}^{m+n} 1\ot e_{ij}\ot t_{ij}(u)(-1)^{\bi\tss\bj+\bj}.
\een
The defining relations for the algebra $\Y(\gl_{m|n})$ take
the form of the $RTT$-{\em relation}
\beql{RTT}
R(u-v)\ts T_1(u)\ts T_2(v)=T_2(v)\ts T_1(u)\ts R(u-v).
\eeq
More explicitly, the defining relations can be written
with the use of super-commutator in terms of the series \eqref{tiju} as
\beql{defrel}
\big[\tss t_{ij}(u),t_{kl}(v)\big]=\frac{1}{u-v}
\big(t_{kj}(u)\ts t_{il}(v)-t_{kj}(v)\ts t_{il}(u)\big)
(-1)^{\bi\tss\bj+\bi\tss\bk+\bj\tss\bk}.
\eeq
It is clear from the defining relations that any
re-labelling map $t_{ij}(u)\mapsto t_{\si(i)\tss\si(j)}(u)$,
preserving the parities of the generators, where
$\si$ is a permutation of the set $\{1,\dots,m+n\}$, extends to
an isomorphism between the respective Yangians. Therefore,
all algebras $\Y(\gl^{\tss\se}_{m|n})$ associated with different
parity sequences $\se$ are isomorphic to each other.

The {\em Yangian} $\Y(\sll_{m|n})$
is the subalgebra of
$\Y(\gl_{m|n})$ which
consists of the elements stable under
the automorphisms
\beql{muf}
t_{ij}(u)\mapsto f(u)\ts t_{ij}(u)
\eeq
for all series
$f(u)\in 1+u^{-1}\CC[[u^{-1}]]$.
If $m\ne n$, then we have the tensor product decomposition
\beql{tensordecom}
\Y(\gl_{m|n})=\ZY(\gl_{m|n})\ot \Y(\sll_{m|n}),
\eeq
where $\ZY(\gl_{m|n})$ denotes the center of $\Y(\gl_{m|n})$. The center is freely
generated by the coefficients of the {\em quantum Berezinian}. Its definition
in the case of standard parity sequence $\se=\se^{\st}$ goes back to \cite{n:qb}; see also
\cite{g:rl} and \cite{g:gd} for more details on the properties of the
Berezinian, including a formula in terms of the Gaussian generators.
Its generalization for an arbitrary parity sequence $\se$ was found in
\cite{t:sa}.

The universal enveloping algebra $\U(\gl_{m|n})$ can be regarded as a subalgebra of
$\Y(\gl_{m|n})$ via the embedding
\beql{emb}
E_{ij}\mapsto t_{ij}^{(1)}(-1)^{\bi},
\eeq
while the mapping
\beql{ev}
t_{ij}(u)\mapsto \de_{ij}+E_{ij}(-1)^{\bi}\ts u^{-1}
\eeq
defines the {\em evaluation homomorphism} $\ev{:}\ts\Y(\gl_{m|n})\to\U(\gl_{m|n})$.

The Yangian $\Y(\gl_{m|n})$ is a Hopf algebra with the coproduct
defined by
\beql{Delta}
\De: t_{ij}(u)\mapsto \sum_{k=1}^{m+n} t_{ik}(u)\ot t_{kj}(u).
\eeq

For any parity sequence $\se$ and a nonnegative integer $p$ introduce
the {\em Yangian of level $p$}, denoted by
$\Y_p(\gl^{\tss\se}_{m|n})$ (or simply by $\Y_p(\gl_{m|n})$ for a fixed $\se$),
as the quotient of the algebra $\Y(\gl^{\tss\se}_{m|n})$ by the two-sided ideal generated by
all elements $t_{ij}^{(r)}$ with $r>p$. For a given $p$,
all algebras $\Y_p(\gl^{\tss\se}_{m|n})$ associated with different
parity sequences $\se$ are isomorphic to each other.

Note that by the results of \cite{bbg:pw} and \cite{p:fw},
the Yangians $\Y_p(\gl_{m|n})$ are isomorphic to
certain finite $\Wc$-algebras corresponding to the general linear Lie superalgebras.

\section{Finite-dimensional irreducible representations}
\label{sec:hw}

A representation $V$ of the algebra $\Y(\gl^{\tss\se}_{m|n})$ associated with a
parity sequence $\se$,
is called a {\em highest weight representation}
if there exists a nonzero vector
$\xi\in V$ such that $V$ is generated by $\xi$,
\begin{alignat}{2}
t_{ij}(u)\ts\xi&=0 \qquad &&\text{for}
\quad 1\leqslant i<j\leqslant m+n, \qquad \text{and}\non\\
t_{ii}(u)\ts\xi&=\la_i(u)\ts\xi \qquad &&\text{for}
\quad i=1,\dots,m+n,
\non
\end{alignat}
for some formal series
\beql{laiu}
\la_i(u)\in 1+u^{-1}\CC[[u^{-1}]].
\eeq
The vector $\xi$ is called the {\em highest vector}
of $V$ and the $(m+n)$-tuple $\la(u)=(\la_{1}(u),\dots,\la_{m+n}(u))$ is called
its {\em highest weight}.

Given an arbitrary tuple $\la(u)=(\la_{1}(u),\dots,\la_{m+n}(u))$
of formal series of the form \eqref{laiu},
the {\em Verma module} $M(\la(u))$ is defined as the quotient of the algebra $\Y(\gl^{\tss\se}_{m|n})$ by
the left ideal generated by all coefficients of the series $t_{ij}(u)$
with $1\leqslant i<j\leqslant m+n$, and $t_{ii}(u)-\la_i(u)$ for
$i=1,\dots,m+n$. We will denote by $L(\la(u))$ its irreducible quotient and will say that
$L(\la(u))$ is associated with the parity sequence $\se$. It is clear that the isomorphism
class of $L(\la(u))$ is determined by $\la(u)$.

\bpr\label{prop:fdhw}
Every finite-dimensional irreducible representation of the algebra $\Y(\gl^{\tss\se}_{m|n})$
is isomorphic to $L(\la(u))$ for a certain highest weight
$\la(u)=(\la_{1}(u),\dots,\la_{m+n}(u))$.
\epr

\bpf
The argument is essentially the same as the proof of the counterpart
of the property for the Yangians associated with Lie algebras as in \cite{zh:sy};
see also \cite[Sec.~3.2]{m:yc} and Remark~\ref{rem:pbwrep}.
\epf

Note a property of the Yangians associated with sequences of length two.
It is easy to check that
we have the isomorphisms
\beql{zeon}
\Y(\gl^{\tss01}_{1|1})\to \Y(\gl^{\tss10}_{1|1}),\qquad t^{\tss01}_{ij}(u)\mapsto t^{\tss10}_{ij}(-u),
\eeq
and
\beql{zeze}
\Y(\gl^{\tss00}_{2|0})\to \Y(\gl^{\tss11}_{0|2}),\qquad t^{\tss00}_{ij}(u)\mapsto t^{\tss11}_{ij}(-u),
\eeq
where the superscripts indicate which algebras
the generating series correspond to.

\bpr\label{prop:necc}
Suppose that the representation $L(\la(u))$ of $\Y(\gl^{\tss\se}_{m|n})$
is finite-dimensional. Then for each subsequence of $\se$
of the form $\se_i\se_{i+1}=01$ or $\se_i\se_{i+1}=10$ the ratio $\la_i(u)/\la_{i+1}(u)$ is the expansion
into a series in $u^{-1}$ of a rational function in $u$,
\ben
\frac{\la_i(u)}{\la_{i+1}(u)}=\frac{\ol Q_i(u)}{Q_i(u)},
\een
for some monic polynomials $\ol Q_i(u)$ and $Q_i(u)$ in $u$ of the same degree.

Moreover,
for each subsequence of $\se$
of the form $\se_i\se_{i+1}=00$ or $\se_i\se_{i+1}=11$
we have
\ben
\frac{\la_i(u)}{\la_{i+1}(u)}=\frac{P_i(u+1)}{P_i(u)}\Fand
\frac{\la_{i+1}(u)}{\la_i(u)}=\frac{P_i(u+1)}{P_i(u)},
\een
respectively, for a monic polynomial $P_i(u)$ in $u$.
\epr

\bpf
Consider the subalgebra $\Y_i$ of $\Y(\gl^{\tss\se}_{m|n})$ generated by the coefficients
of the series $t_{ii}(u)$, $t_{i\ts i+1}(u)$, $t_{i+1\ts i}(u)$ and $t_{i+1,i+1}(u)$.
This subalgebra is isomorphic to one of the four algebras occurring in
\eqref{zeon} and \eqref{zeze}. The cyclic span $\Y_i\tss \xi$ is a finite-dimensional
module with the highest weight $(\la_i(u),\la_{i+1}(u))$.
In view of the isomorphisms \eqref{zeon} and \eqref{zeze}, we only need to derive
the required conditions for the highest weight representations
of the algebras $\Y(\gl^{\tss01}_{1|1})$ and $\Y(\gl^{\tss00}_{2|0})$.
The derivation is the same in both cases following the original approach of \cite{t:im}
as shown in \cite{zh:rs}; see also \cite[Sec.~3.3]{m:yc}.
\epf

The necessary conditions of Proposition~\ref{prop:necc} need not be sufficient in general.
However, they are sufficient in the case
of the standard parity sequence $\se=\se^{\st}$.
The next theorem, whose proof we outline below, is due to Zhang~\cite{zh:sy}. We will write
$\Y(\gl_{m|n})$ for the Yangian associated with the sequence $\se^{\st}=0\dots 0\tss 1\dots 1$.

\bth\label{thm:cri}
The representation $L(\la(u))$ of $\Y(\gl_{m|n})$
is finite-dimensional if and only if there exist monic polynomials
$P_1(u),\dots,P_{m-1}(u),\ol Q_m(u),Q_m(u),P_{m+1}(u),\dots,P_{m+n-1}(u)$ in the variable
$u$ such that
\begin{align}
\frac{\la_i(u)}{\la_{i+1}(u)}&=\frac{P_i(u+1)}{P_i(u)}\qquad\text{for}\quad i=1,\dots,m-1,
\non\\[0.4em]
\frac{\la_{i+1}(u)}{\la_i(u)}&=\frac{P_i(u+1)}{P_i(u)}\qquad\text{for}\quad i=m+1,\dots,m+n-1,
\non\\
\intertext{and}
\frac{\la_m(u)}{\la_{m+1}(u)}&=\frac{\ol Q_m(u)}{Q_m(u)}.
\non
\end{align}
\eth

\bpf
The necessity of the conditions follows from Proposition~\ref{prop:necc}, so we only need to
show that they are sufficient.
The conditions determine the highest weight
$\la(u)$ up to a simultaneous multiplication of all components
$\la_i(u)$ by a series $f(u)\in 1+u^{-1}\CC[[u^{-1}]]$.
This operation corresponds to twisting the action of
the algebra $\Y(\gl_{m|n})$ on $L(\la(u))$ by the automorphism \eqref{muf}.
Therefore, it is enough to prove that
a particular module $L(\la(u))$ corresponding to a given set
of polynomials is finite-dimensional. Moreover, all components $\la_i(u)$
of the highest weight of such a module can be assumed to be polynomials in $u^{-1}$.
Factorize them to write
\ben
\la_i(u)=(1+\la_i^{(1)}\tss(-1)^{\bi}\ts u^{-1})\dots (1+\la_i^{(p)}\tss(-1)^{\bi}\ts u^{-1}),
\qquad i=1,\dots,m+n,
\een
for a certain positive integer $p$.
By the assumptions on the highest weight $\la(u)$, we can renumber the factors
in these decompositions to ensure that all tuples
\ben
\la^{(r)}=(\la_1^{(r)},\dots,\la_{m+n}^{(r)}),\qquad r=1,\dots,p,
\een
satisfy the conditions
\ben
\la_i^{(r)}-\la_{i+1}^{(r)}\in\ZZ_+
\quad\text{for all}\quad i\ne m.
\een
Use the evaluation homomorphism \eqref{ev} to equip
each finite-dimensional irreducible representation $L(\la^{(r)})$
of the Lie superalgebra $\gl_{m|n}$
with an $\Y(\gl_{m|n})$-module
structure.
The coproduct formula \eqref{Delta} implies that
the $\Y(\gl_{m|n})$-module $L(\la(u))$ is isomorphic to a subquotient of the tensor
product of finite-dimensional evaluation modules
\beql{teeva}
L(\la^{(1)})\ot\dots\ot L(\la^{(p)}),
\eeq
thus proving that $L(\la(u))$ is finite-dimensional.
\epf

\bre\label{rem:pbwrep}
The Yangian representation theory relies on a suitable version of the
Poincar\'e--Birkhoff--Witt theorem. Such versions for $\Y(\gl_{m|n})$ were given in \cite{g:gd},
providing a missing part of the proof contained in the earlier work \cite{zh:sy}.
For extensions of the theorem and additional comments on the validity of the previous proofs
for the Drinfeld presentation, see \cite{t:sa}.
\qed
\ere

\section{Odd reflections}
\label{sec:or}

The proof of Theorem~\ref{thm:cri} shows, that after twisting with a suitable automorphism
of the form \eqref{muf}, every finite-dimensional irreducible representation
of the algebra $\Y(\gl_{m|n})$ factors through a representation of
the Yangian $\Y_p(\gl_{m|n})$ of some level $p$. This is implied by the
evaluation and coproduct formulas \eqref{ev} and \eqref{Delta} which show that
all generators $t_{ij}^{(r)}$ of $\Y(\gl_{m|n})$ with $r>p$ act as the zero operators in
the modules \eqref{teeva}. We will now be concerned with
representations of the Yangians $\Y_p(\gl_{m|n})$.

We will use the same notation $t_{ij}^{(r)}$ for the generators of $\Y_p(\gl_{m|n})$ and
extend the definitions of highest weight representations to $\Y_p(\gl_{m|n})$ accordingly.
Introduce the polynomials
\ben
T_{ij}(u)=u^p\ts t_{ij}(u)=\de_{ij}u^p+t_{ij}^{(1)}\tss u^{p-1}+\dots+t_{ij}^{(p)},\qquad i,j=1,\dots,m+n.
\een
The defining relations for $\Y_p(\gl_{m|n})$ in terms of the polynomials $T_{ij}(u)$
take the same form \eqref{defrel}, with the replacement $t_{ij}(u)\mapsto T_{ij}(u)$.

We will now concentrate on the particular case $m=n=1$ and introduce a suitable version of
the odd reflections for the corresponding Yangian of level $p$. We will work with
the Yangian $\Y_p(\gl_{1|1})$ associated with the standard parity sequence
$\se=01$; the results in the opposite case $\se=10$ will then follow
by the application of the isomorphism \eqref{zeon}.
Note the relation
\beql{tco}
(u-v+1)\ts T_{21}(u)\ts T_{21}(v)=-(u-v-1)\ts T_{21}(v)\ts T_{21}(u)
\eeq
implied by \eqref{defrel}. By setting $v=u+1$ we get
\beql{anni}
T_{21}(u+1)\ts T_{21}(u)=0.
\eeq

We will keep the same notation $\la(u)=(\la_1(u),\la_2(u))$ for the highest weight
of the representation $L(\la(u))$ of $\Y_p(\gl_{1|1})$,
but now $\la_1(u)$ and $\la_2(u)$ will be regarded as monic polynomials in $u$ of degree $p$
with $T_{ii}(u)\ts\xi=\la_i(u)\ts\xi$.
To state the next proposition,
write the expansions
\beql{laot}
\la_1(u)=(u+\al_1)\dots (u+\al_p)\Fand \la_2(u)=(u+\be_1)\dots (u+\be_p)
\eeq
for some $\al_i,\be_i\in\CC$.
Assume that $\al_i\ne\be_j$ for all $i,j\in\{1,\dots,p\}$.
Furthermore, to avoid vanishing of some elements below
due to \eqref{anni}, we will impose the following {\em anti-string condition}
on the parameters $\al_1,\dots,\al_p$:
\beql{strco}
\al_i-\al_{j}+1\ne 0 \qquad\text{for all}\quad 1\leqslant i<j\leqslant p.
\eeq
This condition is not restrictive and it can be achieved by a renumbering
of the roots of $\la_1(u)$; see also Remark~\ref{rem:veze}.
Keeping these assumptions, we can now state a key result underlying the use
of odd reflections in the Yangian context.

\bpr\label{prop:prove}
The vector
\beql{ze}
\ze=T_{21}(-\al_1)\dots T_{21}(-\al_p)\ts\xi\in L(\la(u))
\eeq
is nonzero.
Moreover,
the following relations hold:
\begin{align}
T_{11}(u)\ts\ze&=(u+\al_1-1)\dots (u+\al_p-1)\ts\ze,
\label{too}\\[0.4em]
T_{22}(u)\ts\ze&=(u+\be_1-1)\dots (u+\be_p-1)\ts\ze,
\label{ttt}\\
\intertext{and}
T_{21}(u)\ts\ze&=0.
\label{tto}
\end{align}
\epr

\bpf
We have
\ben
T_{12}(u)\ts\ze=\sum_{i=1}^p(-1)^{i-1}\ts
T_{21}(-\al_1)\dots \big[T_{12}(u),T_{21}(-\al_i)\big]
\dots T_{21}(-\al_p)\ts\xi.
\een
By the defining relations,
\ben
\big[T_{12}(u),T_{21}(-\al_i)\big]=-\frac{1}{u+\al_i}
\big(T_{22}(u)\ts T_{11}(-\al_i)-T_{22}(-\al_i)\ts T_{11}(u)\big)
\een
and
\beql{tooa}
T_{11}(u)\ts  T_{21}(-\al_j)=\frac{u+\al_j-1}{u+\al_j}\ts T_{21}(-\al_j)\ts T_{11}(u)
+\frac{1}{u+\al_j}\ts T_{21}(u)\ts T_{11}(-\al_j).
\eeq
Since
\beql{ttoxi}
T_{11}(u)\ts\xi=(u+\al_1)\dots (u+\al_p)\ts\xi,
\eeq
an easy induction yields
\ben
T_{12}(u)\ts\ze=\sum_{i=1}^p(-1)^{i}\ts\prod_{l=1}^{i-1}(u+\al_l)\ts
\prod_{l=i+1}^{p}(u+\al_l-1)\ts
T_{21}(-\al_1)\dots T_{22}(-\al_i)
\dots T_{21}(-\al_p)\ts\xi.
\een
Now set $u=-\al_p+1$ to get
\ben
T_{12}(-\al_p+1)\ts\ze=(-1)^p\ts \prod_{i=1}^p(\be_i-\al_p)\ts\prod_{l=1}^{p-1}(\al_l-\al_p+1)\ts
T_{21}(-\al_1)\dots T_{21}(-\al_{p-1})\ts \xi.
\een
The assumptions on the parameters and the anti-string condition \eqref{strco}
imply that the numerical coefficient of the vector on the right hand side is nonzero.
Therefore, by applying the operators $T_{12}(-\al_{p-1}+1),\dots,T_{12}(-\al_1+1)$
to the vector on the right hand side consecutively,
we obtain the highest vector $\xi$ with a nonzero coefficient, thus proving
that the vector $\ze$ is nonzero.

To prove the second part of the proposition, note that relation
\eqref{too} is implied by \eqref{tooa} and \eqref{ttoxi},
as was pointed out in the above
calculation. Furthermore, relation \eqref{tto} is a consequence
of the identity in $\Y_p(\gl_{1|1})$ which holds for independent
variables $u_1,\dots,u_{p+1}$:
\beql{idep}
T_{21}(u_1)\dots T_{21}(u_{p+1})=0.
\eeq
Indeed, the expression on the left hand side is a polynomial in these variables with
the degree with respect to each variable not exceeding $p-1$. However,
the expression vanishes at the evaluations $u_i=u_j+1$
for all $1\leqslant i<j\leqslant p+1$ due to \eqref{tco} and \eqref{anni}, and therefore must be
identically zero. Alternatively, an equivalent form of \eqref{idep} reads
$
t_{21}^{(r_1)}\dots t_{21}^{(r_{p+1})}=0
$
for all $r_i\in\{1,\dots,p\}$, which holds by \eqref{tco} as the generators $t_{21}^{(r)}$
and $t_{21}^{(s)}$ anticommute modulo lower degree terms.

Finally, for the proof of \eqref{ttt} we use the quantum Berezinian for the Yangian $\Y(\gl_{1|1})$,
as introduced in \cite{n:qb}; see also \cite{g:rl}.
We will regard
$L(\la(u))$ as the $\Y(\gl_{1|1})$-module, obtained via
the composition with the natural epimorphism $\Y(\gl_{1|1})\to \Y_p(\gl_{1|1})$.
All coefficients of the series $b(u)$ defined by
\ben
b(u)=\big(t_{22}(u)-t_{21}(u)\ts t_{11}(u)^{-1}t_{12}(u)\big)\ts t_{11}(u)^{-1}
\een
belong to the center of the algebra $\Y(\gl_{1|1})$.
They act as multiplications
by scalars in the
module $L(\la(u))$ which are found by $b(u)\mapsto \la_2(u)/\la_1(u)$.
On the other hand, the defining relations imply that
\ben
b(u)\ts t_{11}(u)\ts t_{11}(u+1)=t_{22}(u+1)\ts t_{11}(u)+t_{12}(u+1)\ts t_{21}(u).
\een
Apply the series on both sides
to the vector $\ze\in L(\la(u))$
by using the already
verified relations \eqref{too} and \eqref{tto}
to derive \eqref{ttt}.
\epf

\bre\label{rem:veze}
Note that since the generators $t_{21}^{(r)}$
and $t_{21}^{(s)}$ anticommute modulo lower degree terms, the vector $\ze$ can be written
as $\ze=c\ts t_{21}^{(1)}\dots t_{21}^{(p)}\ts\xi$ for a nonzero constant $c$.
This also implies that $\ze$ must be an eigenvector for $T_{11}(u)$ and $T_{22}(u)$, although
expression \eqref{ze} appears to be more convenient for the calculation of the eigenvalues.
\qed
\ere

We will also state a generalization of Proposition~\ref{prop:prove} to the case where
the polynomials $\la_1(u)$ and $\la_2(u)$ have common roots.
Consider the expansions \eqref{laot} and
renumber the parameters if necessary, to assume that for some $k\in\{0,1,\dots,p\}$
we have $\al_i=\be_i$ for $i=k+1,\dots,p$, whereas $\al_i\ne\be_j$ for all $i,j\in\{1,\dots,k\}$.
The anti-string condition \eqref{strco} will now be assumed for
the parameters $\al_1,\dots,\al_k$.

\bco\label{cor:prove}
The fraction
\beql{kfra}
\ol T_{21}(u)=\frac{T_{21}(u)}{\ga(u)}\qquad\text{with}\quad
\ga(u)=(u+\al_{k+1})\dots (u+\al_p)
\eeq
is a polynomial in $u$, as an operator in $L(\la(u))$.
Moreover, the vector
\ben
\ze=\ol T_{21}(-\al_1)\dots \ol T_{21}(-\al_k)\ts\xi\in L(\la(u))
\een
is nonzero and satisfies $T_{21}(u)\ts\ze=0$,
\begin{align}
T_{11}(u)\ts\ze&=(u+\al_1-1)\dots (u+\al_k-1)\ts\ga(u)\ts\ze,
\non\\[0.4em]
T_{22}(u)\ts\ze&=(u+\be_1-1)\dots (u+\be_k-1)\ts\ga(u)\ts\ze.
\non
\end{align}
\eco

\bpf
Consider the representation $L(\bar\la(u))$ of the Yangian $\Y_k(\gl_{1|1})$ of level $k$, where the
components of the highest weight have the form
\ben
\bar\la_1(u)=(u+\al_1)\dots (u+\al_k)\Fand \bar\la_2(u)=(u+\be_1)\dots (u+\be_k).
\een
Denote the generator polynomials for the Yangian $\Y_k(\gl_{1|1})$ by $\ol T_{ij}(u)$
and use the same symbols for the corresponding operators on $L(\bar\la(u))$.
It is immediate from the defining relations \eqref{defrel} that the assignment
\ben
T_{ij}(u)\mapsto\ga(u)\ts \ol T_{ij}(u)\qquad\text{for all}\quad i,j\in\{1,2\}
\een
defines a representation of the Yangian $\Y_p(\gl_{1|1})$
on the vector space $L(\bar\la(u))$. This representation is clearly irreducible
and isomorphic to $L(\la(u))$. This proves the first part of the proposition,
because the action of the operator in $L(\la(u))$ defined in \eqref{kfra}
corresponds to the action of the operator
$\ol T_{21}(u)$ in $L(\bar\la(u))$. The remaining parts of the corollary are now
immediate from Proposition~\ref{prop:prove}.
\epf

Now consider the Yangian $\Y_p(\gl^{\tss\se}_{m|n})$ associated with
an arbitrary parity sequence $\se$ and introduce some notation to state the main results.
Let $L(\la(u))$ be the irreducible highest weight representation
of $\Y_p(\gl^{\tss\se}_{m|n})$ with a certain highest weight $\la(u)$, whose components,
regarded as polynomials in $u$ with $T_{jj}(u)\ts\xi=\la_j(u)\ts\xi$, are
given by
\ben
\la_j(u)=(u+\la_j^{(1)})\dots (u+\la_j^{(p)}),\qquad j=1,\dots,m+n,\quad \la_j^{(p)}\in\CC.
\een

Suppose that the parity sequence $\se$ has a subsequence of
the form $\se_i\se_{i+1}=10$. Then we let $\se^+$ denote the sequence obtained from $\se$
by replacing this subsequence with $01$ and leaving the remaining terms unchanged.
Similarly, if $\se$ has a subsequence
$\se_i\se_{i+1}=01$, we let $\se^-$ denote the sequence obtained from $\se$
by replacing this subsequence with $10$ and leaving the remaining terms unchanged.

With a chosen value of $i$, we will assume that the roots of the polynomials
$\la_i(u)$ and $\la_{i+1}(u)$ are numbered in such a way that
\beql{laeq}
\la_i^{(r)}=\la_{i+1}^{(r)}\qquad \text{for all}\quad r=k+1,\dots,p
\eeq
for certain $k\in\{0,1,\dots,p\}$, while $\la_i^{(r)}\ne\la_{i+1}^{(s)}$
for all $1\leqslant r,s\leqslant k$. Then define new polynomials $\la^{\pm}_j(u)$
for $j=1,\dots,m+n$ by
\ben
\bal
\la^{\pm}_i(u)&=(u+\la_{i+1}^{(1)}\pm 1)\dots (u+\la_{i+1}^{(k)}\pm 1)
(u+\la_{i+1}^{(k+1)})\dots (u+\la_{i+1}^{(p)}),\\[0.4em]
\la^{\pm}_{i+1}(u)&=(u+\la_i^{(1)}\pm 1)\dots (u+\la_i^{(k)}\pm 1)
(u+\la_i^{(k+1)})\dots (u+\la_i^{(p)}),
\eal
\een
and $\la^{\pm}_j(u)=\la_j(u)$ for $j\ne i,i+1$. These polynomials form the respective
highest weights $\la^+(u)$ and $\la^-(u)$ for the Yangian of level $p$.

\bth\label{thm:oddre}
With the above assumptions, the following holds.
\begin{enumerate}
\item
If the parity sequence $\se$ has a subsequence
$\se_i\se_{i+1}=01$
then the $\Y_p(\gl_{m|n})$-module $L(\la(u))$ associated with $\se$ is isomorphic to
the $\Y_p(\gl_{m|n})$-module $L(\la^-(u))$ associated with $\se^-$.
\item
If the parity sequence $\se$ has a subsequence
$\se_i\se_{i+1}=10$
then the $\Y_p(\gl_{m|n})$-module $L(\la(u))$ associated with $\se$ is isomorphic to
the $\Y_p(\gl_{m|n})$-module $L(\la^+(u))$ associated with $\se^+$.
\end{enumerate}
\eth

\bpf
The two parts of the theorem are equivalent, and so it is sufficient to prove the first part.
The subalgebra $\Y_i$ of $\Y_p(\gl_{m|n})$ generated by the coefficients
of the polynomials $T_{ii}(u)$, $T_{i\ts i+1}(u)$, $T_{i+1\ts i}(u)$ and $T_{i+1,i+1}(u)$
is isomorphic to $\Y_p(\gl_{1|1})$. The cyclic span $\Y_i\tss \xi$ of the highest vector
$\xi$ of $L(\la(u))$ is a highest weight module
associated with the parity sequence $01$, with the highest weight
$(\la_i(u),\la_{i+1}(u))$.
Renumber the parameters $\la_i^{(1)},\dots,\la_i^{(k)}$,
if necessary,
to satisfy the anti-string condition \eqref{strco}.
Apply Corollary~\ref{cor:prove} to get the nonzero vector
\ben
\ze_i=\ol T_{i+1\ts i}(-\la_i^{(1)})\dots \ol T_{i+1\ts i}(-\la_i^{(k)})\ts\xi\in L(\la(u)),
\een
where
\ben
\ol T_{i+1\ts i}(u)=\frac{T_{i+1\ts i}(u)}{\ga_i(u)}\qquad\text{with}\quad
\ga_i(u)=(u+\la_i^{(k+1)})\dots (u+\la_i^{(p)}).
\een
An easy induction with the use of
the defining relations \eqref{defrel} shows that for the vector $\ze_i$ we have
\beql{relab}
T_{ab}(u)\ts \ze_i=0
\eeq
for all $a<b$ with $a<i$ or $b>i+1$.
Twist the action of $\Y_p(\gl^{\tss\se}_{m|n})$
on $L(\la(u))$ with the
re-labelling isomorphism $T_{ab}(u)\mapsto T_{\si_i(a)\tss\si_i(b)}(u)$
for the transposition $\si_i=(i,i+1)$
to get a representation of $\Y_p(\gl^{\tss\se^-}_{m|n})$
on the same vector space.
Corollary~\ref{cor:prove} and relations \eqref{relab} imply
that the resulting $\Y_p(\gl^{\tss\se^-}_{m|n})$-module
is an irreducible highest weight representation with the highest vector $\ze_i$
and the highest weight $\la^-(u)$, so it
is isomorphic to $L(\la^-(u))$.
\epf

Due to the homomorphisms \eqref{emb} and \eqref{ev},
the Yangian $\Y_1(\gl_{m|n})$ of level $1$ is isomorphic to the universal
enveloping algebra $\U(\gl_{m|n})$. Therefore Theorem~\ref{thm:oddre} for $p=1$
is a particular case of the well-known properties of odd reflections; see e.g.
\cite{ps:cg}.

Any two parity sequences are related by a chain of transitions $01\leftrightarrow 10$
involving two adjacent entries. Therefore, as
a consequence of Theorem~\ref{thm:oddre} and the criterion of
Theorem~\ref{thm:cri}, we can get necessary and sufficient conditions for
the irreducible highest weight
representation $L(\la(u))$ of $\Y_p(\gl_{m|n})$ associated with
an arbitrary parity sequence $\se$ to be finite-dimensional. As with the case $p=1$,
an explicit form of such conditions should look quite complicated; cf.
\cite[Sec.~2.4]{cw:dr}, where extremal weights in polynomial modules are discussed.

\bex\label{ex:glot}
Take the parity sequence $\se=101$ for $\Y_p(\gl_{1|2})$ and consider the
irreducible highest weight representation $L(\la(u))$ with
\ben
\bal
\la_1(u)&=(u+\al_1)\dots (u+\al_p),\\[0.3em]
\la_2(u)&=(u+\be_1)\dots (u+\be_p),\\[0.3em]
\la_3(u)&=(u+\ga_1)\dots (u+\ga_p).
\eal
\een
Renumber the roots of the polynomials $\la_1(u)$ and $\la_2(u)$ if necessary,
and introduce the parameter $k\in\{0,1,\dots,p\}$ to have
$\al_i=\be_i$ for $i=k+1,\dots,p$ and $\al_i\ne\be_j$ for all $i,j\in\{1,\dots,k\}$.
By Part 2 of Theorem~\ref{thm:oddre}, the representation $L(\la(u))$
is isomorphic to the $\Y_p(\gl_{1|2})$-module $L(\la^+(u))$ associated with the standard parity
sequence $\se^+=011$, where
\ben
\bal
\la^+_1(u)&=(u+\be_1+1)\dots (u+\be_k+1)(u+\be_{k+1})\dots (u+\be_p),\\[0.3em]
\la^+_2(u)&=(u+\al_1+1)\dots (u+\al_k+1)(u+\al_{k+1})\dots (u+\al_p),\\[0.3em]
\la^+_3(u)&=(u+\ga_1)\dots (u+\ga_p).
\eal
\een
Hence, by Theorem~\ref{thm:cri}, the representation $L(\la(u))$ is finite-dimensional
if and only if
\ben
\frac{\la^+_{3}(u)}{\la^+_2(u)}=\frac{P(u+1)}{P(u)}
\een
for a monic polynomial $P(u)$ in $u$.
\qed
\eex

Since the odd reflections of Theorem~\ref{thm:oddre} only depend on one subsequence $\se_i\se_{i+1}$,
we can derive a slightly more general transition rule which applies to representations of the Yangian
$\Y(\gl^{\tss\se}_{m|n})$. Consider the irreducible highest weight
$\Y(\gl^{\tss\se}_{m|n})$-module $L(\la(u))$,
where the components of the highest weight are series of the form \eqref{laiu}.
Suppose that the parity sequence $\se$ has a subsequence
$\se_i\se_{i+1}=01$ or $\se_i\se_{i+1}=10$, and that the components $\la_i(u)$ and $\la_{i+1}(u)$
are polynomials in $u^{-1}$,
\ben
\bal
\la_i(u)&=(1+\la_i^{(1)} u^{-1})\dots (1+\la_i^{(k)} u^{-1}),\\[0.3em]
\la_{i+1}(u)&=(1+\la_{i+1}^{(1)} u^{-1})\dots (1+\la_{i+1}^{(k)} u^{-1}),
\eal
\een
without common roots. Transform the highest weight $\la(u)$ by setting
\ben
\bal
\la^{\pm}_i(u)&=\big(1+(\la_{i+1}^{(1)}\pm1)\tss u^{-1}\big)\dots
\big(1+(\la_{i+1}^{(k)}\pm1)\tss  u^{-1}\big),\\[0.3em]
\la^{\pm}_{i+1}(u)&=\big(1+(\la_i^{(1)}\pm1\big)\tss  u^{-1})\dots
\big(1+(\la_i^{(k)}\pm1)\tss  u^{-1}\big),
\eal
\een
and $\la^{\pm}_j(u)=\la_j(u)$ for $j\ne i,i+1$, and keep the notation
$\se^-$ and $\se^+$ used in Theorem~\ref{thm:oddre}.

\bco\label{cor:oddrege}
The $\Y(\gl_{m|n})$-module $L(\la(u))$ associated with the parity sequence
$\se$ is isomorphic to
the $\Y(\gl_{m|n})$-module $L(\la^{\pm}(u))$ associated with $\se^{\pm}$,
where the plus and minus signs are chosen for $\se_i\se_{i+1}=10$ and $\se_i\se_{i+1}=01$,
respectively.
\eco

\bpf
Consider the $\Y(\gl_{1|1})$-subalgebra
of $\Y(\gl_{m|n})$ generated by the coefficients
of the series $t_{ii}(u)$, $t_{i\ts i+1}(u)$, $t_{i+1\ts i}(u)$ and $t_{i+1,i+1}(u)$.
The cyclic span $\ol L=\Y(\gl_{1|1})\tss \xi$ is
a highest weight representation of
$\Y(\gl_{1|1})$ with the highest weight $(\la_i(u),\la_{i+1}(u))$.
By the assumptions, both $\la_i(u)$ and $\la_{i+1}(u)$ are polynomials in $u^{-1}$
of degree not exceeding $k$. By the results of
\cite{zh:rs} describing the structure of such representations,
all vectors $t_{ab}^{(r)}\ts\xi$
with $a,b\in\{i,i+1\}$ and $r>k$ are zero in $L(\la(u))$,
and the representation $\ol L$ factors through the
representation of the Yangian
$\Y_k(\gl_{1|1})$ of level $k$. The corollary then follows by
the application of the relations provided in
Proposition~\ref{prop:prove}, as in the proof of Theorem~\ref{thm:oddre}.
\epf

To conclude, we note that the transition rule of Corollary~\ref{cor:oddrege} is applicable
in a more general case, where
the ratio $\la_i(u)/\la_{i+1}(u)$ is a rational function.
By twisting the module $L(\la(u))$ with a suitable automorphism \eqref{muf}, we can ensure
that the components $f(u)\la_i(u)$ and
$f(u)\la_{i+1}(u)$ of the new highest weight satisfy the required assumptions.
Therefore, taking into account the necessary conditions
of Proposition~\ref{prop:necc}, we can state the general inductive rule to determine whether the given
$\Y(\gl^{\tss\se}_{m|n})$-module $L(\la(u))$ is finite-dimensional:

\begin{itemize}
\item[(i)] if the ratio $\la_i(u)/\la_{i+1}(u)$ is not a rational function for
a certain subsequence
$\se_i\se_{i+1}=10$ or $\se_i\se_{i+1}=01$, then $L(\la(u))$ is infinite-dimensional;
otherwise, go to step (ii);
\item[(ii)] if $\se$ is standard, go to step (iii); otherwise
apply the transition rule of Corollary~\ref{cor:oddrege}
for a subsequence $\se_i\se_{i+1}=10$, then set $\se:=\se^+$ and $\la(u):=\la^+(u)$
and repeat step (i);
\item[(iii)] apply Theorem~\ref{thm:cri}.
\end{itemize}

\bigskip\bigskip

\small

\noindent
School of Mathematics and Statistics\newline
University of Sydney,
NSW 2006, Australia\newline
alexander.molev@sydney.edu.au

\end{document}